\documentclass[10pt]{article}
\usepackage{amsmath,amsthm,amssymb}
\usepackage{graphicx}
\usepackage{color}

\newcommand{\arcsinh}{{\mathrm{arcsinh}}}

\title{Inequalities for Generalized Trigonometric and Hyperbolic Sine Functions}
\author{\small Miao-Kun Wang$^{1}$, Yu-Ming Chu$^{2,*}$ and Yue-Ping Jiang$^3$}
\date{}

\begin{document}
\maketitle
\renewcommand{\thefootnote}{\fnsymbol{footnote}}
{\footnotesize\rm
\noindent $^1$Department of Mathematics, Huzhou Teachers College, Huzhou 313000, China;\\
$^2$School of Mathematics and Computation Sciences, Hunan City University, Yiyang 413000, China;\\
$^3$College of Mathematics and Econometrics, Hunan University,
Changsha
410082, China.\\
Correspondence should be addressed to Yu-Ming Chu,
chuyuming@hutc.zj.cn}

\medskip
\noindent{\bf Abstract}:  We prove that the inequalities
$\sin_{p,q}(\sqrt{rs})\geq \sqrt{\sin_{p,q}(r)\sin_{p,q}(s)}$ and
$\sinh_{p,q}(\sqrt{r^*s^*}) \leq
\sqrt{\sinh_{p,q}(r^*)\sinh_{p,q}(s^*)}$ hold for all
$p,q\in(1,\infty)$, $r,s\in(0,\int_{0}^{1}(1-t^q)^{-1/p}dt)$ and
$r^*,s^*\in(0,\int_{0}^{\infty}(1+t^q)^{-1/p}dt)$, where
$\sin_{p,q}$ and $\sinh_{p,q}$ are the generalized trigonometric and
hyperbolic sine functions, respectively. As a consequence of the
results, we prove a conjecture due to Bhayo and Vuorinen [J. Approx.
Theory, 164(2012)].



\bigskip
\bigskip
\noindent{\bf\large 1. Introduction}
\bigskip
\setcounter{section}{1} \setcounter{equation}{0}

It is well known from basic calculus that
\begin{equation*}
\frac{\pi}{2}=\int_{0}^{1}\frac{1}{\sqrt{1-t^2}}dx
\end{equation*}
and
\begin{equation*}
\arcsin(x)=\int_{0}^{x}\frac{1}{\sqrt{1-t^2}}dt,\quad 0\leq x\leq 1.
\end{equation*}
Since the function $\arcsin(x)$ is a differentiable function on
$[0,1]$ and $t\rightarrow 1/\sqrt{1-t^2}$ is strictly increasing on
$[0,1)$, we can define $\sin$ on $[0,\pi/2]$ as the inverse function
of $\arcsin$. By standard extension procedures we can define the sin
function on $(-\infty,\infty)$.

For $p,q>1$, let
\begin{equation*}
F_{p,q}(x)=\int_{0}^{x}(1-t^q)^{-1/p}dt,\quad x\in[0,1],
\end{equation*}
\begin{equation*}
\frac{\pi_{p,q}}{2}=\int_{0}^{1}(1-t^q)^{-1/p}dt.
\end{equation*}
Then $F_{p,q}:[0,1]\rightarrow [0,\pi_{p,q}/2]$ is an increasing
homeomorphism, denoted by $\arcsin_{p,q}$. Thus its inverse
\begin{equation*}
\sin_{p,q}=F_{p,q}^{-1}
\end{equation*}
is defined on the interval $[0,\pi_{p,q}/2]$. By the similar
extension as the sine function, we can get a differentiable function
$\sin_{p,q}$ defined on $\mathbb{R}$. We call $\sin_{p,q}$ the
generalized $(p,q)$-trigonometric sine function.

We also defined  $\arccos_{p,q}(x)=\arcsin_{p,q}((1-x^p)^{1/q})$
(See [2, 7]), and the inverse of generalized $(p,q)$-hyperbolic sine
function
\begin{equation*}
\arcsinh_{p,q}(x)=\int_{0}^{x}(1+t^q)^{-1/p}dt,\quad x\in(0,\infty).
\end{equation*}

Their inverse functions are
\begin{align*}
&\sin_{p,q}:(0,\pi_{p,q}/2)\rightarrow (0,1),\qquad &\cos_{p,q}:(0,\pi_{p,q}/2)\rightarrow (0,1),\nonumber\\
&\sinh_{p,q}:(0,m^*_{p,q})\rightarrow (0,\infty),\qquad
&m^*_{p,q}=\int_{0}^{\infty}(1+t^q)^{-1/p}dt.
\end{align*}

When $p=q$, the $(p,q)$-functions $\sin_{p,q}$, $\cos_{p,q}$,
$\sinh_{p,q}$, $\arcsin_{p,q}$, $\arccos_{p,q}$ and $\arcsinh_{p,q}$
reduce to $p$-functions $\sin_{p}$, $\cos_{p}$, $\sinh_{p}$,
$\arcsin_{p}$, $\arccos_{p}$ and $\arcsinh_{p}$ (See [5, 8, 10]),
respectively. In particular, when $p=q=2$, the $(p,q)$-functions
become our familiar trigonometric and hyperbolic functions.

Recently, the generalized trigonometric and hyperbolic functions
($(p,q)$-functions and $p$-functions) have been found many important
applications in differential equations, the theory of operator,
approximation theory and other related fields [4, 9, 11].

In face of the importance of generalized trigonometric and
hyperbolic functions they have been studied by many authors from
different points of view [1, 3-8, 10, 12]. Ednumds, Gurka and Lang
[7] gave the basic properties of generalized $(p,q)$-trigonometric
functions, and proved that
\begin{equation*}
\sin_{4/3,4}(2x)=\frac{2\sin_{4/3,4}(x)(\cos_{4/3,4}(x))^{1/3}}{(1+4(\sin_{4/3,4}(x))^4(\cos_{4/3,4}(x))^{4/3})^{1/2}}.
\end{equation*}
Kl\'{e}n, Vuorinen and Zhang [8] generalized  some classical
inequalities for trigonometric and hyperbolic functions, such as
Mitrinovi\'{c}-Adamovi\'{c} inequality and Lazarevi\'{c}'s
inequality.

Bhayo and Vuorinen [2] found that the functions $\arcsin_{p,q}$ and
$\arcsinh_{p,q}$ can be expressed in terms of Gaussian
hypergeometric functions. Applying the vast available information
about the hypergeometric functions, some remarkable properties and
inequalities for generalized trigonometric and hyperbolic functions
are obtained. Moreover, they raised the following conjecture.

\medskip
{\bf Conjecture 1.1.} If $p,q\in(1,\infty)$ and $r,s\in(0,1)$, then
\begin{equation*}
\sin_{p,q}(\sqrt{rs})\geq \sqrt{\sin_{p,q}(r)\sin_{p,q}(s)},
\end{equation*}
\begin{equation*}
\sinh_{p,q}(\sqrt{rs})\leq \sqrt{\sinh_{p,q}(r)\sinh_{p,q}(s)}.
\end{equation*}

The main purpose of this paper is to give a positive answer to the
Conjecture 1.1 and generalized the inequalities in Conjecture 1.1.
Our main result is the following Theorem 1.1.

\medskip
{\bf Theorem 1.1.} If $p,q\in(1,\infty)$, then

(1) Inequality
\begin{equation}
\sin_{p,q}(\sqrt{rs}) \geq \sqrt{\sin_{p,q}(r)\sin_{p,q}(s)}
\end{equation}
holds for all $r,s\in(0,{\pi}_{p,q}/2)$.

(2) Inequality
\begin{equation}
\sinh_{p,q}(\sqrt{r^*s^*}) \leq
\sqrt{\sinh_{p,q}(r^*)\sinh_{p,q}(s^*)}
\end{equation}
holds for all $r^*,s^*\in(0,m_{p,q}^*)$.\\

\bigskip
\bigskip
\noindent{\bf\large 2. Proof of Theorem 1.1}
\bigskip
\setcounter{section}{2} \setcounter{equation}{0}

In order to prove Theorem 1.1 and Conjecture 1.1, we present three
Lemmas at first.

\medskip
{\bf Lemma 2.1.} If $p,q\in(1,\infty)$, then inequality
\begin{equation}
\arcsin_{p,q}(x)>\frac{px(1-x^q)^{1-1/p}}{(q-p)x^q+p}
\end{equation}
holds for $x\in (0,1)$.

\medskip
{\em Proof.} Let
\begin{equation}
\zeta(x)=\arcsin_{p,q}(x)-\frac{px(1-x^q)^{1-1/q}}{(q-p)x^q+p},\quad
x\in(0,1).
\end{equation}
Then simple computations lead to
\begin{equation}
\zeta(0)=0,
\end{equation}
\begin{align}
\zeta'(x)=&\frac{1}{(1-x^q)^{1/p}}-\frac{p^2+p(2q-2p-q^2)x^q+(q-p)^2x^{2q}}{(1-x^q)^{1/p}(p-px^q+q{x^q})^2}\nonumber\\
=&\frac{pq^2x^q}{(1-x^q)^{1/p}(p-px^q+qx^q)^2}>0
\end{align}
for $x\in (0,1)$.

Therefore, Lemma 2.1 follows easily from (2.2)-(2.4).

\medskip
{\bf Lemma 2.2.} If $p,q\in(1,\infty)$, then inequality
\begin{equation}
\frac{x}{\arcsinh_{p,q}(x)}>\frac{(p-q)x^q+p}{p(1+x^q)^{1-1/p}}
\end{equation}
holds for $x\in (0,\infty)$.

\medskip
{\em Proof.} We divide the proof into two cases.

{\bf Case 1} $p\geq q$. Let
\begin{equation}
\eta(x)=\arcsinh_{p,q}(x)-\frac{px(1+x^q)^{1-1/p}}{(p-q)x^q+p},\quad
x\in(0,\infty),
\end{equation}
then simple computations lead to
\begin{equation}
\eta(0)=0,
\end{equation}
\begin{align}
\eta'(x)=&\frac{1}{(1+x^q)^{1/p}}-\frac{p^2+p(2p-2q+q^2)x^q+(q-p)^2x^{2q}}{(1+x^q)^{1/p}(p-qx^q+px^q)^2}\nonumber\\
=&-\frac{pq^2x^q}{(1+x^q)^{1/p}(p-qx^q+px^q)^2}<0
\end{align}
for $x\in (0,\infty)$.

Therefore, inequality (2.5) follows easily from (2.6)-(2.8) and
$p\geq q$.

\medskip
{\bf Case 2} $p<q$. Let $x_{0}=[p/(q-p)]^{1/q}$, then
$(p-q)x^{q}+p>0$ for $x\in(0,x_{0})$ and $(p-q)x^{q}+p<0$ for
$x\in(x_{0},\infty)$. Thus it is sufficient to prove inequality
(2.5) for $x\in(0,x_{0})$, which easily follows from the proof of
Case 1.

\medskip
{\bf Lemma 2.3.} If $p,q>1$, then
$m_{p,q}^*=\int\limits_{0}^{\infty}(1+t^q)^{-1/p}dt>1$.

\medskip
{\em Proof.} From the basic properties of generalized integrals  we
clearly see that $m_{p,q}^*=+\infty$ if $p\geq q$ and
$m_{p,q}^*<+\infty$ if $p<q$. Since
\begin{equation*}
(1+t^{(q-2p)/p})^p>1+t^{q-2p}>1+t^q
\end{equation*}
for $t\in(0,1)$, we get
\begin{align*}
m_{p,q}^*=&\int_{0}^{\infty}\frac{1}{(1+t^q)^{1/p}}dt=\int_{0}^{1}\frac{1}{(1+t^q)^{1/p}}dt+\int_{1}^{\infty}\frac{1}{(1+t^q)^{1/p}}dt\\
=&\int_{0}^{1}\frac{1}{(1+t^q)^{1/p}}dt+\int_{0}^{1}\frac{t^{q/p-2}}{(1+t^q)^{1/p}}dt=\int_{0}^{1}\frac{1+t^{q/p-2}}{(1+t^q)^{1/p}}dt>1.
\end{align*}

\bigskip
{\bf\em Proof of Theorem 1.1.} For part (1), without loss of
generality, we assume that $0<x\leq y<1$. Define
\begin{equation}
J(x,y)=\frac{{\arcsin_{p,q}(H_{p_{1}}(x,y))}^2}{\arcsin_{p,q}(x)\arcsin_{p,q}(y)},\quad
p_{1}\in \mathbb{R},
\end{equation}
where
\begin{equation*}
H_{p^*}(a,b)=\left\{
\begin{array}{ll}
\left(\frac{a^{p^*}+b^{p^*}}{2}\right)^{{1}/{p^*}},&p^*\neq0,\\
\sqrt{ab},&p^*=0
\end{array}\right.
\end{equation*}
is the H\"{o}lder mean of order $p^*\in\mathbb{R}$ of two positive
numbers $a$ and $b$.

Let $t=H_{p_{1}}(x,y)$, then ${\partial t}/{\partial
x}=\left({x}/{t}\right)^{p_{1}-1}/2$. If $x<y$, then $t>x$. By
logarithmic differentiation, we get
\begin{equation}
\frac{1}{J(x,y)}\frac{\partial J}{\partial
x}=x^{p_{1}-1}\left[\frac{t^{1-p_{1}}}{\arcsin_{p,q}(t)(1-t^q)^{1/p}}-\frac{x^{1-p_{1}}}{\arcsin_{p,q}(x)(1-x^q)^{1/p}}\right].
\end{equation}

Let
\begin{equation}
F(x)=\frac{x^{1-p_{1}}}{\arcsin_{p,q}(x)(1-x^q)^{1/p}},\quad
x\in(0,1),
\end{equation}
then logarithmic differentiation yields
\begin{align}
\frac{F'(x)}{F(x)}=&(1-p_{1})\frac{1}{x}-\frac{1}{\arcsin_{p,q}(x)(1-x^q)^{1/p}}+\frac{qx^{q-1}}{p(1-x^q)}\nonumber\\
=&\frac{1}{x}(G(x)+1-p_{1}),
\end{align}
where
\begin{equation*}
G(x)=\frac{qx^{q}}{p(1-x^q)}-\frac{x}{\arcsin_{p,q}(x)(1-x^q)^{1/p}}.
\end{equation*}
Note that
\begin{equation}
\lim_{x\rightarrow 0}G(x)=-1, \quad \lim_{x\rightarrow
1}G(x)=+\infty.
\end{equation}

It follows from Lemma 2.1 and (2.13) that the range of $G(x)$ is
$(-1,\infty)$. Thus from (2.12) we conclude that $F'(x)\geq 0$ for
$x\in(0,1)$ if and only if $p_{1}\leq 0$. Namely, $F(x)$ is strictly
increasing on $(0,1)$ if and only if $p_{1}\leq 0$. Moreover, if
$p_{1}>0$, then $F(x)$ is not monotone.

Next, we divide the proof into two cases.

{\bf Case A} $p_{1}\leq 0$. Then from equation (2.10) and the
monotonicity of $F(x)$ we clearly see that $\partial J/\partial
x>0$. Hence $J(x,y)\leq J(y,y)=1$. Then by (2.9) we get
\begin{equation}
\arcsin_{p,q}(H_{p_{1}}(x,y))\leq
\sqrt{\arcsin_{p,q}(x)\arcsin_{p,q}(y)}
\end{equation}
for $p_{1}\leq 0$, with equality if and only if $x=y$.

{\bf Case B} $p_{1}>0$. Then using the similar argument in Case A,
we conclude that there exists $x_{1},x_{2},y_{1},y_{2}\in(0,1)$ such
that
\begin{equation*}
\arcsin_{p,q}(H_{p_{1}}(x_{1},y_{1}))<
\sqrt{\arcsin_{p,q}(x_{1})\arcsin_{p,q}(y_{1})},
\end{equation*}
\begin{equation*}
\arcsin_{p,q}(H_{p_{1}}(x_{2},y_{2}))>
\sqrt{\arcsin_{p,q}(x_{2})\arcsin_{p,q}(y_{2})}.
\end{equation*}

Finally, taking $x=\sin_{p,q}(r)$ and $y=\sin_{p,q}(s)$ in
inequality (2.14), we have
\begin{equation}
\sin_{p,q}(\sqrt{rs})\geq H_{p_{1}}(\sin_{p,q}(r),\sin_{p,q}(s))
\end{equation}
for all $r,s\in(0,\pi_{p,q}/2)$ if and only if  $p_{1}\leq 0$. In
particular, inequality (1.1) follows from (2.15) with $p_{1}=0$.\\

For part (2), without loss of generality, we assume that $0<x\leq
y<\infty$. Define
\begin{equation}
J^*(x,y)=\frac{{\arcsinh_{p,q}(H_{p_{2}}(x,y))}^2}{\arcsinh_{p,q}(x)\arcsinh_{p,q}(y)},\quad
p_{2}\in \mathbb{R}.
\end{equation}
Let $t=H_{p_{2}}(x,y)$, then ${\partial t}/{\partial
x}=\left({x}/{t}\right)^{p_{2}-1}/2$. If $x<y$, then $t>x$. By
logarithmic differentiation, we get
\begin{equation}
\frac{1}{J^*(x,y)}\frac{\partial J^*}{\partial
x}=x^{p_{2}-1}\left[\frac{t^{1-p_{2}}}{\arcsinh_{p,q}(t)(1+t^q)^{1/p}}-\frac{x^{1-p_{2}}}{\arcsinh_{p,q}(x)(1+x^q)^{1/p}}\right].
\end{equation}

Let
\begin{equation}
F^*(x)=\frac{x^{1-p_{2}}}{\arcsin_{p,q}(x)(1+x^q)^{1/p}},\quad
x\in(0,\infty),
\end{equation}
then logarithmic differentiation leads to
\begin{align}
\frac{{F^*}'(x)}{F^*(x)}=&(1-p_{2})\frac{1}{x}-\frac{1}{\arcsinh_{p,q}(x)(1+x^q)^{1/p}}-\frac{qx^{q-1}}{p(1+x^q)}\nonumber\\
=&\frac{1}{x}(1-p_{2}-G^*(x)),
\end{align}
where
\begin{equation*}
G^*(x)=\frac{x}{\arcsinh_{p,q}(x)(1+x^q)^{1/p}}+\frac{qx^{q}}{p(1+x^q)}.
\end{equation*}
Note that
\begin{equation}
\lim_{x\rightarrow 0}G^*(x)=1.
\end{equation}

It follows from Lemma 2.2 that $G^*(x)>1$ for $x\in(0,\infty)$. Then
we have $\inf\limits_{x\in(0,\infty)}G^{*}(x)=1$. Hence equation
(2.19) leads to the conclusion  that ${F^*}'(x)\leq 0$ for
$x\in(0,\infty)$ if and only if $p_{2}\geq 0$. Namely, $F^*(x)$ is
strictly decreasing on $(0,\infty)$ if and only if $p_{2}\geq 0$.
Hence if $p_{2}\geq 0$, then $\partial J^{*}/\partial x<0$ by
(2.17), and $J^*(x,y)\geq J^*(y,y)=1$. Then from (2.16) we have
\begin{equation}
\arcsinh_{p,q}(H_{p_{2}}(x,y))\geq
\sqrt{\arcsinh_{p,q}(x)\arcsin_{p,q}(y)}
\end{equation}
for $p_{2}\geq 0$, with equality if and only if $x=y$.

Taking $x=\sinh_{p,q}(r^*)$ and $y=\sinh_{p,q}(s^*)$ in inequality
(2.21), we get
\begin{equation}
\sinh_{p,q}(\sqrt{r^*s^*})\leq
H_{p_{2}}(\sinh_{p,q}(r^*),\sinh_{p,q}(s^*))
\end{equation}
for all $r^*,s^*\in(0,m_{p,q})$ if $p_{2}\geq 0$. In particular,
inequality (1.2) follows from (2.22) with $p_{2}=0$.

{\bf Remark 2.1.} Conjecture 1.1  follows easily from
$\pi_{p,q}/2>1$ and Lemma 2.3 together with Theorem 1.1.

\medskip
\noindent{\bf\normalsize Acknowledgements}

This research was supported by the Natural Science Foundation of
China (Grant Nos. 11071059, 11071069, 11171307), and the Innovation
Team Foundation of the Department of Education of Zhejiang Province
(Grant No. T200924).

\medskip
\def\refname

\end{document}